\input amstex
\documentstyle{amsppt}
\magnification 1000

\def\diam{\operatorname{diam}}
\def\R{\Bbb R}

\def\N{\Bbb N}
\def\Z{\Bbb Z}

\def\A{\Cal A}
\def\U{\Cal U}
\def\V{\Cal V}
\def\Ord{\operatorname{Ord}}

\def\dim{\operatorname{dim}}

\def\asdim{\operatorname{asdim}}

\def\trasdim{\operatorname{trasdim}}
\def\coasdim{\operatorname{coasdim}}
\NoBlackBoxes

\topmatter
\title
A metric space with  transfinite  asymptotic dimension $2\omega$
\endtitle

\author
T.Radul
\endauthor
\address
Kazimierz Wielki University, Bydgoszcz (Poland) and Ivan Franko National University of Lviv (Ukraine)
 \newline
e-mail: tarasradul\@yahoo.co.uk
\endaddress

\abstract We build an example of a metric  space with  transfinite  asymptotic dimension $2\omega$.
\endabstract
\endtopmatter
\document
\baselineskip18pt

{\bf 0. Introduction.} The asymptotic dimension $\asdim$ of a metric space was
defined by Gromov for studying asymptotic invariants of discrete
groups [1]. This dimension can be considered as an asymptotic
analogue of the Lebesgue covering dimension $\dim$.

A transfinite extensions $\trasdim$ for the
asymptotic dimensions  $\asdim$ was introduced in [2]. There was  constructed a proper metric space $X$ such that $\trasdim
X=\omega$. Metric spaces  $X_{\omega+k}$ with $\trasdim
X_{\omega+k}=\omega+k$ were constructed in [3] for each natural $k$.

{\bf 1. Preliminaries.}

 A family $\A$ of
subsets of a metric space is called {\it uniformly bounded} if
there exists a number $C>0$ such that $\diam A\le C$ for each
$A\in\A$; $\A$ is called $r$-{\it disjoint} for some $r>0$ if
$d(A_1,A_2)\ge r$ for each $A_1$, $A_2\in\A$ such that $A_1\neq
A_2$.

The {\it asymptotic dimension} of a metric space $X$ does not
exceed $n\in\N\cup\{0\}$ (written $\asdim X\le n$) iff for every
$D>0$ there exists a uniformly bounded cover $\U$ of $X$ such that
$\U=\U_0\cup\dots\cup\U_n$, where all $\U_i$ are $D$-disjoint [1].

Since the definition of $\asdim$ is not inductive, we cannot
immediately extend this dimension to transfinite numbers. We need some set-theoretical
construction used by Borst to extend the covering dimension and metric
dimension [4,5].

Let $L$ be an arbitrary set. By $Fin L$ we shall denote the
collection of all finite, non-empty subsets of $L$. Let $M$ be a
subset of $Fin L$. For $\sigma\in\{\emptyset\}\cup Fin L$ we put
$$M^\sigma=\{\tau\in Fin L\mid \sigma\cup\tau\in M\ \text{and}\
\sigma\cap\tau=\emptyset\}.$$ Let $M^a$ abbreviate $M^{\{a\}}$ for
$a\in L$.

 Define the ordinal number $\Ord M$ inductively
as follows

$\Ord M=0$ iff $M=\emptyset$,

$\Ord M\le\alpha$ iff for every $a\in L$, $\Ord M^a<\alpha$,

$\Ord M=\alpha$ iff $\Ord M\le\alpha$ and $\Ord M<\alpha$ is not
true, and

$\Ord M=\infty$ iff $\Ord M>\alpha$ for every ordinal number
$\alpha$.

Given a metric space $(X,d)$, let us define the following collection:
$$A(X,d)=\{\sigma\in Fin \N\mid \ \text{ there is no uniformly bounded  families}\
\V_i\ \text{for}\ i\in\sigma$$ $$\text{such that}\
\cup_{i\in\sigma}\V_i\ \text{covers}\ X\ \text{and}\ \V_i\
\text{is}\ i-\text{disjoint}\}.$$

Let $(X,d)$ be a metric space. Then put $\trasdim X=\Ord A(X,d)$ [2].

We will need some properties of $\trasdim$ proved in  [2].

\proclaim {Proposition A,[2]} Let $X$ be a metric space and $Y\subset
X$. Then $\trasdim Y\le\trasdim X$.
\endproclaim

We denote by $N_R(A)=\{x\in X|d(x,A)\le R\}$ for a metric space
$X$, $A\subset X$ and $R>0$.

\proclaim {Proposition B,[2]} Let $Y$ be a metric space and $X\subset Y$.
Then $\trasdim N_n(X)=\trasdim X$ for each $n\in\N$.
\endproclaim

We also will use another transfinite extension $\coasdim$ of $\asdim$ defined in [3].

We represent every ordinal number $\gamma=\lambda(\gamma)+n(\gamma)$ where $\lambda(\gamma)$ is the limit ordinal number or $0$ and $n(\gamma)\in \N\cup\{0\}$. For a metric space $X$ we define complementary-finite asymptotic dimension $\coasdim(X)$ inductively
as follows [3]:

$\coasdim(X)=-1$ iff $X=\emptyset$,

$\coasdim(X)\le\gamma$ iff for every $r>0$ there exist $r$-disjoint uniformly bounded families $\V_0,\dots,\V_{n(\gamma)}$ of subsets of $X$ such that $\coasdim(X\bigcup(\bigcup_{i=0}^{n(\gamma)}))<\lambda(\gamma)$,

$\coasdim(X)=\gamma$ iff $\coasdim(X)\le\gamma$ and $\coasdim(X)\le\beta$ is not
true for each $\beta<\gamma$, and

$\coasdim(X)=\infty$ iff $\coasdim(X)\le\alpha$ is not true for every ordinal number
$\alpha$.

It was proved in [3] that inequality  $\coasdim(X)\le\omega+k$ implies $\trasdim(X)\le\omega+k$ for each metric space. By routine checking we can generalize this fact.

\proclaim {Proposition 1} We have $\trasdim(X)\le\coasdim(X)$ for each metric space $X$.
\endproclaim

Using similar arguments as in the proof of Proposition B, we can prove the following statement.

\proclaim {Proposition 2} Let $Y$ be a metric space and $X\subset Y$.
Then $\coasdim N_n(X)=\coasdim X$ for each $n\in\N$.
\endproclaim

{\bf 3. The main result.}

We are going to construct  a  metric space
$X_{2\omega}$ such that $\trasdim(X_{2\omega})=\coasdim(X_{2\omega})=2\omega$. We develop the ideas from [2] and [3]. Firstly let us recall the construction of the spaces $X_{\omega+k}$ for $k\in\N$ from [3].

 We consider $\R^n$ with the sup-metric defined as follows
$d_\infty((k_1,\dots,k_n),(l_1,\dots,l_n))=\max\{|k_1-l_1|,\dots,|k_n-l_n|\}$. For $k\le n$ we consider the natural isometrical embedding $i_n^k:\R^k\to\R^n$ defined by the formula  $i_n^k(x_1,\dots,x_k)=(x_1,\dots,x_k,0,\dots,0)$.

For $i,k\in\N$ we define $$X_k^{(i)}=\{x\in\R^i\mid|\{l\in\{1,\dots,i\}\mid x_l\notin 2^i\Z\}|\le k\}.$$
Put $X_{\omega+k}=\bigsqcup_{i=1}^\infty X_k^{(i)}$. Define a metric $d_k$ on $X_{\omega+k}$ as follows. Consider any $x$, $y\in X_{\omega+k}$. Let $x\in X_k^{(m)}$ and $y\in X_k^{(n)}$ with $m\le n$. Put 
$c_k=0$ if $n=m$ and $c_k=k(m+(m+1)+\dots+(n-1))$ if $m<n$. Define $d_k(x,y)=\max\{d_\infty(i_n^m(x),y),c_k\}$. It is proved in [3] that $\trasdim(X_{\omega+k})=\coasdim(X_{\omega+k})=\omega+k$ for each $k\in\N$. Let us remark that we slightly modified the metric $d_k$ compared with [3], but this modification does not affect on the statements from [3].

It is easy to see that for each pair $n$, $m\in \N$ such that  $m\le n$ the natural including $j_n^m:X_{\omega+m}\to X_{\omega+n}$ is an isometrical embedding.

Put $Y_k^{(m)}=X_k^{(m)}\cap (2^k\Z)^m$ and $Y_{\omega+k}=\bigsqcup_{i=1}^\infty Y_k^{(i)}$. It is easy to see that $N_{2^k}(Y_{\omega+k})=X_{\omega+k}$, hence $\trasdim(Y_{\omega+k})=\omega+k$ for each $k\in\N$. Put $X_{2\omega}=\bigsqcup_{i=1}^\infty Y_{\omega+i}$. Define a metric $d$ on $X_{2\omega}$ as follows. Consider any $x$, $y\in X_{2\omega}$. Let $x\in X_k^{(m)}$ and $y\in X_k^{(n)}$ with $m\le n$. Put
$c=0$ if $n=m$ and $c=m+(m+1)+\dots+(n-1)$ if $m<n$. Define $d(x,y)=\max\{d_n(j_n^m(x),y),c\}$.

\proclaim {Theorem 1} $\trasdim(X_{2\omega})=2\omega$.
\endproclaim

\demo{Proof} The inequality $\trasdim(X_{2\omega})\ge 2\omega$ follows
from the fact that $\trasdim(Y_{\omega+k})=\omega+k$ for each $k\in\N$ and Proposition A.

Let us show that $\coasdim(X_{2\omega})\le2\omega$. Let $r>0$. Choose $k\in\N$ such that $k\ge r$. Consider the family $\V=\{\{x\}|x\in\cup_{i=k}^\infty Y_{\omega+i}\}$. The family $\V$ is evidently $r$-disjoint. We have $X_{2\omega}\setminus \cup\V=\cup_{i=1}^{k-1} Y_{\omega+i}$. But $\cup_{i=1}^{k-1} Y_{\omega+i}\subset N_c(Y_{\omega+k-1})$, where $c=1+\dots+k-1$. Hence $\coasdim(X_{2\omega}\setminus \cup\V)\le\omega+k-1<2\omega$.

Now, the theorem follows from Proposition 1.
\enddemo

Certainly, the presented result is only one step ahead. The general problem still remains: for each countable ordinal $\xi$   find a metric space $X_\xi$ with $\trasdim(X_\xi)=\xi$.


\Refs \baselineskip12pt
\parskip3pt

\ref \no 1 \by M.Gromov \book Asymptotic invariants of infinite
groups. Geometric group theory. v.2 \publaddr Cambridge University
Press \yr 1993
\endref

\ref \no 2 \by T.Radul \paper On transfinite extension of asymptotic dimension \jour Topology Appl. \yr 2010 \vol 157 \pages 2292–2296
\endref

\ref \no 3 \by Yan Wu, Jingming Zhu \paper A metric space with its transfinite asymptotic
dimension $\omega+k$  \yr 2019 \jour  arXiv:1912.02103v1
\endref

\ref \no 4 \by P.Borst \paper Classification of weakly
infinite-dimensional spaces  \jour Fund. Math \yr 1988 \vol 130
\pages 1-25 \endref

\ref \no 5 \by P.Borst \paper Some remarks concerning $C$-spaces
\jour Topology Appl. \yr 2007 \vol 154
\pages 665-674 \endref








\endRefs
\enddocument